\documentclass[11pt]{amsart}
\usepackage{latexsym}
\usepackage{amssymb}
\usepackage{amsfonts}
\usepackage{graphicx}
\usepackage{epsf}

\newtheorem{theorem}{Theorem}[section]
\newtheorem{corollary}[theorem]{Corollary}

\newtheorem{proposition}[theorem]{Proposition}

\newtheorem{remark}[theorem]{Remark}

\def\bea{\begin{eqnarray*}}
\def\eea{\end{eqnarray*}}

\def\ot{\otimes}
\def\ra{\rightarrow}

\def\bea{\begin{eqnarray*}}
\def\eea{\end{eqnarray*}}

\begin{document}
\title{A linear algebra approach to graded Frobenius algebras }

\author{
S. D\u{a}sc\u{a}lescu$^{1}$, C. N\u{a}st\u{a}sescu$^{2}$, L.
N\u{a}st\u{a}sescu$^{2}$ and P. Rebenciuc$^{1}$}
\address{$^1$ University of Bucharest, Faculty of Mathematics and Computer Science,
Str. Academiei 14, Bucharest 1, RO-010014, Romania} \address{ $^2$
Institute of Mathematics of the Romanian Academy, PO-Box
1-764\\
RO-014700, Bucharest, Romania}
\address{
 e-mail: sdascal@fmi.unibuc.ro, Constantin\_nastasescu@yahoo.com,
 lauranastasescu@gmail.com, paul.rebenciuc@s.unibuc.ro
}

\date{}
\maketitle

\begin{abstract}
If $A$ is a finite-dimensional algebra graded by a group $G$, and
$\sigma \in G$, we define a variant of paratrophic matrix associated
with $A$ and $\sigma$, and we use it to characterize the
$\sigma$-graded Frobenius property for $A$. We discuss the
invertibility of such paratrophic matrices, and then use them to
check whether certain graded algebras are $\sigma$-graded Frobenius
or (graded) symmetric. As an application, we uncover (graded)
Frobenius and symmetric properties of Koszul duals of quantum
polynomial algebras. We derive a structure result for
$\sigma$-graded Frobenius algebras by only using linear algebra
methods. \\[1mm]
2020 MSC: 16W50, 15A09, 16L60, 16S36, 16S37\\[1mm]
Key words: graded algebra, graded Frobenius algebra, graded
symmetric algebra, paratrophic matrix, quantum polynomial algebra,
Koszul dual, good grading on a matrix algebra.
\end{abstract}

\section{Introduction and preliminaries}

Frobenius algebras form one of the most interesting classes of
finite-dimensional algebras. Originated in the work of Frobenius
\cite{Frobenius}, these objects have occurred in algebra, topology,
knot theory, Hopf algebra theory, mathematical physics, etc. A first
characterization of these algebras was given by Frobenius in terms
of certain matrices called paratrophic matrices.

A finite-dimensional algebra $A$ is called Frobenius if $A$ and its
linear dual $A^*$ are isomorphic as left $A$-modules. A version of
the Frobenius property was introduced for graded algebras in
\cite{dnn1}: if $A$ is a finite-dimensional algebra graded by a
group $G$, and $\sigma \in G$, then $A$ is called $\sigma$-graded
Frobenius if $A^*\simeq A(\sigma)$ as graded left $A$-modules, where
$A(\sigma)$ is the shift of $A$ by $\sigma$; see the definitions in
Section \ref{sectiongradedFrobenius}. Such objects have occurred  in
classification results of certain algebras in non-commutative
geometry \cite{lpwz}, and for proving a non-commutative
Bernstein-Gelfand-Gelfand correspondence \cite{jor}. There is a
special interest in $\varepsilon$-graded Frobenius algebras, where
$\varepsilon$ is the neutral element of $G$; these are simply called
graded Frobenius algebras, and they are just the Frobenius algebras
in the monoidal category of $G$-graded vector spaces, see
\cite{dnn1}. If the strongest property that $A^*\simeq A$ as graded
$A$-bimodules is satisfied, then $A$ is called graded symmetric.

Our aim is to develop a linear algebra approach to graded Frobenius
algebras in the spirit of Frobenius' initial approach. Let $G$ be a
multiplicative group with neutral element $\varepsilon$, and let
$A=\oplus_{g\in G}A_g$ be a finite-dimensional $G$-graded algebra
over a field $K$. For each $g\in G$, let $(e_i)_{i\in J_g}$ be a
linear basis of $A_g$; if $A_g=0$, we take $J_g=\emptyset$. Denote
by $J$ the union of all $J_g$ over all $g\in G$; then $(e_i)_{i\in
J}$ is a homogeneous basis of $A$. For any $i,j\in J$, write
$e_ie_j=\sum_{\ell \in J}c_{ij\ell}e_{\ell}$; the scalars
$(c_{ij\ell})_{i,j,\ell\in J}$ are called the structure constants of
$A$.

Fix some $\sigma\in G$. For any family of scalars
$\alpha=(\alpha_\ell)_{\ell\in J_\sigma}$, we consider the matrix
$P(\sigma,\alpha)\in M_J(K)$ whose $(i,j)$-entry is $\sum_{\ell\in
J_\sigma}\alpha_\ell c_{ij\ell}$ for any $i,j\in J$; we call such a
$P(\sigma,\alpha)$ a $\sigma$-paratrophic matrix associated with
$A$. We use the convention that a sum indexed by the empty set is
$0$, thus a $\sigma$-paratrophic matrix is automatically $0$ if
$A_\sigma=0$.

We also consider for each $g\in G$ the matrix $C_g$ with rows
indexed by $J_{\sigma g^{-1}}\times J_\sigma$ and columns indexed by
$J_g$, whose entry on the $(r,\ell)$th row and $j$th column is
$c_{rj\ell}$, for any $r\in J_{\sigma g^{-1}}$, $\ell \in J_\sigma$
and $j\in J_g$. By convention, a matrix with an empty set of rows or
columns is an empty matrix, is invertible and has rank 0. We denote
by $P(\sigma,\alpha)_{F,F'}$ the submatrix of $P(\sigma,\alpha)$
obtained by intersecting rows in $F$ with columns in $F'$. We give
in Section \ref{sectioninvertibleparatrophic} the following
characterization of the invertibility
of a $\sigma$-paratrophic matrix.\\

{\bf Theorem A.} {\it Let $\alpha=(\alpha_\ell)_{\ell\in J_\sigma}$
be a family of scalars. The following are equivalent.\\
{\rm (1)} $P(\sigma,\alpha)$ is invertible.\\
{\rm (2)} $|J_\varepsilon|=|J_\sigma|$,
$P(\sigma,\alpha)_{J_\varepsilon,J_\sigma}$ is invertible, and ${\rm
rank}\, C_g=|J_g|$ for any $g\in G$.\\
{\rm (3)} For any $g\in G$ we have $|J_{\sigma g^{-1}}|=|J_g|$ and
$P(\sigma,\alpha)_{J_{\sigma g^{-1}},J_g}$ is invertible.}\\

Paratrophic matrices are used in Section
\ref{sectiongradedFrobenius} to prove the following characterization
of the $\sigma$-graded Frobenius property.\\

{\bf Proposition B.} {\it $A$ is $\sigma$-graded Frobenius if and
only if there exists a family of scalars
$\alpha=(\alpha_\ell)_{\ell\in J_\sigma}$ such that the
$\sigma$-paratrophic matrix
$P(\sigma,\alpha)$ is invertible.}\\

We also give a version of Proposition B for graded symmetric
algebras. We use these results in Section \ref{sectiuneexemple} for
testing the $\sigma$-graded Frobenius and the (graded) symmetric
properties on a class of algebras constructed as follows. Let $n\geq
2$ be an integer and let ${\bf q}=(q_{ij})_{1\leq i<j\leq n}$ be a
family of non-zero scalars. Consider the quantum polynomial algebra
$K_{\bf q}[X_1,\ldots,X_n]$ in $n$ indeterminates; this is the
quotient of the free $K$-algebra $K\langle X_1,\ldots,X_n\rangle$ in
$n$ letters by the ideal generated by all $X_jX_i-q_{ij}X_iX_j$,
with $1\leq i<j\leq n$. Then consider the finite-dimensional algebra
$A({\bf q})=K_{\bf q}[X_1,\ldots,X_n]/(X_1^2,\ldots,X_n^2)$. If
$q_{ij}=-1$ for any $i,j$, then $A({\bf q})$ is isomorphic to the
exterior algebra of a vector space of dimension $n$. We explain that
$A({\bf q})$ is the Koszul dual of a certain algebra of quantum
polynomials.

Denoting by $x_1,\ldots,x_n$ the classes of $X_1,\ldots,X_n$ in
$A({\bf q})$, a basis of $A({\bf q})$ consists of all $
x_1^{\varepsilon_1}\ldots x_n^{\varepsilon_n}$ with
$\varepsilon_1,\ldots,\varepsilon_n\in \{ 0,1\}$. Then $A({\bf q})$
has a $\mathbb{Z}$-grading whose homogeneous component of degree $d$
is the span of all $x_1^{\varepsilon_1}\ldots x_n^{\varepsilon_n}$
with $\varepsilon_1+\ldots +\varepsilon_n=d$ for any $d$. There is
also a $\mathbb{Z}_2^n$-grading on $A({\bf q})$, such that the
homogeneous component of degree
$(\widehat{\varepsilon_1},\ldots,\widehat{\varepsilon_n})$ is
$Kx_1^{\varepsilon_1}\ldots x_n^{\varepsilon_n}$ for any
$\varepsilon_1,\ldots,\varepsilon_n\in \{ 0,1\}$.
We prove the following.\\

{\bf Proposition C.} {\it {\rm (i)} $A({\bf q})$ is a Frobenius
algebra, and it is a symmetric algebra if and only if
$$q_{in}=\left( \prod_{1\leq p<i}q_{pi}\right)\cdot \left( \prod_{i<p\leq
n-1}q_{ip}^{-1}\right)$$ for any $1\leq i\leq n-1$.\\
{\rm (ii)} Regarded with the $\mathbb{Z}_2^n$-grading, $A({\bf q})$
is $(\hat{1},\ldots,\hat{1})$-graded Frobenius, but it is not
$\sigma$-graded Frobenius for any $\sigma\neq
(\hat{1},\ldots,\hat{1})$.\\
{\rm (iii)} Regarded with the $\mathbb{Z}$-grading, $A({\bf q})$ is
$n$-graded Frobenius, but it is not $j$-graded Frobenius for any
$j\neq n$.}\\

We also use the symmetric version of Proposition B to show that a
full matrix algebra endowed with a good grading, i.e., a grading
such that any usual matrix unit is a homogeneous element, is graded
symmetric.

Finally, in Section \ref{sectionstructureresult} we derive from
Theorem A the following structure result for
$\sigma$-graded Frobenius algebras.\\

{\bf Corollary D.} (\cite[Theorem 4.3]{dnn1}) {\it A
finite-dimensional $G$-graded algebra $A$ is $\sigma$-graded
Frobenius if and only if it is left $\sigma$-faithful and
$(A_\sigma)^*\simeq A_\varepsilon$ as
left $A_\varepsilon$-modules.}\\

The original proof of this result, given in \cite{dnn1}, uses
constructions and methods of graded rings, including the coinduced
functor. Our approach only uses linear algebra techniques.

For basic facts on graded algebras we refer \cite{nvo}.

\section{Invertible $\sigma$-paratrophic matrices}
\label{sectioninvertibleparatrophic}

We keep the setup described in the Introduction. Thus $G$ is a group
with neutral element $\varepsilon$, $\sigma$ is a fixed element of
$G$, $A$ is a finite-dimensional $G$-graded algebra, $(e_i)_{i\in
J}$ is a basis of $A$ consisting of homogeneous elements, with the
ones of degree $g\in G$ being indexed by $J_g$, and
$(c_{ij\ell})_{i,j,\ell}$ are the associated structure constants. If
$i\in J_g$ and $j\in J_h$, where $g,h\in G$, then $e_ie_j$ has
degree $gh$, so $c_{ij\ell}=0$ for any $\ell \notin J_{gh}$.

If $g,h,k\in G$ and $i\in J_g, j\in J_h, p\in J_k$, we obtain from
$(e_ie_j)e_p=e_i(e_je_p)$ that

\begin{equation} \label{equation*}
\sum_{\ell\in J_{gh}}c_{ij\ell}c_{\ell pr}=\sum_{\ell\in
J_{hk}}c_{i\ell r}c_{jp\ell}
\end{equation}
for any $r\in J$.\\

 {\bf Proof of Theorem A:} For any $i,j\in J$,
denote by $p_{ij}=\sum_{\ell \in J_\sigma}\alpha_\ell c_{ij\ell}$
the $(i,j)$th entry of $P(\sigma,\alpha)$.

${\rm (1)}\Rightarrow {\rm (2)}$ The partition $J=\cup_{g\in G}J_g$
induces a presentation of $P(\sigma,\alpha)$ into possibly non-zero
blocks $P(\sigma,\alpha)_{J_{\sigma g^{-1}},J_g}$, and other blocks
of the form $P(\sigma,\alpha)_{J_h,J_g}$, all of them zero. If
$|J_\varepsilon|\neq |J_\sigma|$, then ${\rm
det}(P(\sigma,\alpha))=0$ by using Laplace expansion over rows $i$
with $i\in J_\varepsilon$ in the case where $|J_\varepsilon|>
|J_\sigma|$, respectively over columns $j$ with $j\in J_\sigma$, in
the case where $|J_\varepsilon|< |J_\sigma|$. Thus $|J_\varepsilon|=
|J_\sigma|$, and Laplace expansion over rows in $J_\varepsilon$
shows that ${\rm det}(P(\sigma,\alpha))=M{\rm
det}(P(\sigma,\alpha)_{J_\varepsilon,J_\sigma})$ for some $M\in K$,
so ${\rm det}(P(\sigma,\alpha)_{J_\varepsilon,J_\sigma})\neq 0$,
showing that $P(\sigma,\alpha)_{J_\varepsilon,J_\sigma}$ is
invertible.

Pick now $g\in G$ and let $(y_j)_{j\in J_g}\subset K$ such that

\begin{equation}\label{equation****}
\sum_{j\in J_g}c_{rj\ell}y_j=0
\end{equation}
for any $(r,\ell)\in J_{\sigma g^{-1}}\times J_\sigma$. Fix some
$r\in J_{\sigma g^{-1}}$. Then

\bea \sum_{j\in J_g}p_{rj}y_j&=&\sum_{j\in J_g}\sum_{\ell\in
J_\sigma}\alpha_\ell c_{rj\ell} y_j\\
&=&\sum_{\ell \in J_\sigma}\left( \alpha_\ell \sum_{j\in
J_g}c_{rj\ell}y_j\right)\\
&=&0 \;\;\;\;\;\;\;\;\;\;\;\; (\mbox{by equation
(\ref{equation****})}) \eea Now defining $y_j=0$ for any $j\in
J\setminus J_g$, we obtain that $(y_j)_{j\in J}$ is a solution of
the homogeneous system whose matrix is $P(\sigma,\alpha)$. Since
this matrix is invertible, we get that $y_j=0$ for each $j\in J$. We
conclude that the system consisting of the equations
(\ref{equation****}), with $(r,\ell)\in J_{\sigma g^{-1}}\times
J_\sigma$, has just the trivial solution, and so ${\rm rank}\,
C_g=|J_g|$.

${\rm (2)}\Rightarrow {\rm (3)}$ Let $g\in G$, and let $(y_j)_{j\in
J_g}\subset K$ be such

\begin{equation}\label{equationsist3}\sum_{j\in J_g}p_{rj}y_j=0
\end{equation}
for any $r\in J_{\sigma g^{-1}}$. Thus

$$\sum_{j\in J_g}\sum_{\ell\in J_\sigma}\alpha_\ell c_{rj\ell}
y_j=0$$ for each $r\in J_{\sigma g^{-1}}$. Then clearly

$$\sum_{r\in J_{\sigma g^{-1}}}c_{qr'r}\left(\sum_{j\in J_g}\sum_{\ell \in
J_\sigma}\alpha_\ell c_{rj\ell} y_j\right)=0$$ for any $r'\in
J_{\sigma g^{-1}}$ and $q\in J_\varepsilon$. Reorganizing the sums,
we get

\begin{equation} \label{equation**}
\sum_{\ell \in J_\sigma}\sum_{j\in J_g}\left(\sum_{r\in J_{\sigma
g^{-1}}}c_{qr'r}c_{rj\ell}\right)\alpha_\ell y_j=0
\end{equation}

By (\ref{equation*}), we have

$$\sum_{r\in J_{\sigma
g^{-1}}}c_{qr'r}c_{rj\ell}=\sum_{\ell '\in
J_\sigma}c_{q\ell'\ell}c_{r'j\ell'},$$ and replacing in
(\ref{equation**}), we obtain
$$\sum_{\ell \in J_\sigma}\sum_{j\in J_g}\sum_{\ell'\in
J_\sigma}c_{q\ell'\ell}c_{r'j\ell'}\alpha_\ell y_j=0.$$ This means
that
$$\sum_{\ell'\in J_\sigma}p_{q\ell'}\left(\sum_{j\in
J_g}c_{r'j\ell'}y_j\right)=0$$ for any $q\in J_\varepsilon$. As
$P(\sigma,\alpha)_{J_\varepsilon,J_\sigma}$ is invertible, we obtain
that
$$\sum_{j\in
J_g}c_{r'j\ell'}y_j=0$$ for any $\ell'\in J_\sigma$ and $r'\in
J_{\sigma g^{-1}}$. This is a homogeneous system in $(y_j)_{j\in
J_g}$, whose matrix is $C_g$. Since ${\rm rank}\, C_g=|J_g|$, the
system has just the trivial solution. Hence the system with
equations (\ref{equationsist3}) has only trivial solution, showing
that ${\rm rank}\, P(\sigma,\alpha)_{J_{\sigma g^{-1}},J_g}=|J_g|$.
In particular
\begin{equation}\label{equation***}
|J_g|\leq |J_{\sigma g^{-1}}|
\end{equation}
for any $g\in G$. Summing over all $g\in G$, we find $\sum_{g\in
G}|J_g|\leq \sum_{g\in G}|J_{\sigma g^{-1}}|$. As both sides of the
inequality are equal to $|J|$, we must have equality in
(\ref{equation***}) for each $g$, and then clearly,
$P(\sigma,\alpha)_{J_{\sigma g^{-1}},J_g}$ is invertible.

${\rm (3)}\Rightarrow {\rm (1)}$ $P(\sigma,\alpha)$ is invertible,
since it consists of invertible square blocks
$P(\sigma,\alpha)_{J_{\sigma g^{-1}},J_g}$, where $(J_g)_{g\in G}$
and $(J_{\sigma g^{-1}})_{g\in G}$ are partitions of $J$. \qed

\section{$\sigma$-graded Frobenius algebras and paratrophic
matrices} \label{sectiongradedFrobenius}

We will need the following result which describes the situations
where the linear dual of a right module over a finite-dimensional
algebra is isomorphic to the regular left representation of the
algebra. In the particular case where the module is just the regular
right representation of the algebra, this is just Frobenius'
characterization of Frobenius algebras. If $U$ is a right module
over a $K$-algebra $R$, then the linear dual $U^*$ of $U$ is a left
$R$-module with the left action, denoted by $\rightharpoonup$,
induced by the left $R$-action on $U$.

\begin{proposition} \label{propUdual}
Let $R$ be a $K$-algebra and $U$ a right $R$-module such that ${\rm
dim}\, R={\rm dim}\, U=m$. Let $\varepsilon_1,\ldots,\varepsilon_m$
be a basis of $R$, $u_1,\ldots,u_m$ be a basis of $U$, and let
$(c_{ij\ell})_{1\leq i,j,\ell\leq m}\subset K$ such that
\begin{equation}\label{structureU}
u_i\varepsilon_j=\sum_{1\leq \ell\leq m}c_{ij\ell}u_\ell
\end{equation}
for any $1\leq i,j\leq m$. Then $U^*\simeq R$ as left $R$-modules if
and only if there exists a family $(\alpha_\ell)_{1\leq \ell\leq
m}\subset K$ satisfying the property
$$(*)\;\;\;\;\;\;\;\mbox{ the matrix } \left( \sum_{1\leq \ell\leq
m}\alpha_\ell c_{ij\ell}\right)_{1\leq i,j\leq m}\;\;\mbox{ is
invertible}$$ Moreover, in this case, any isomorphism $\varphi:R\ra
U^*$ is of the form $\varphi(r)=r\rightharpoonup u^*$, for some
$u^*\in U^*$ such that the family $(u^*(u_\ell))_{1\leq \ell\leq m}$
satisfies the property $\rm (*)$.
\end{proposition}
\begin{proof}
Any morphism of left $R$-modules $\varphi:R\ra U^*$ is of the form
$\varphi(r)=r\rightharpoonup u^*$ for any $r\in R$, for some $u^*\in
U^*$. Looking at dimensions, such a $\varphi$ is an isomorphism if
and only if it is injective.

Let $r=\sum_{1\leq j\leq m}y_j\varepsilon_j$. Then $r \in {\rm
Ker}\, \varphi$ if and only if $(r\rightharpoonup u^*)(u_i)=0$ for
any $1\leq i\leq m$, which rewrites as $u^*(\sum_{1\leq j\leq
m}y_ju_i\varepsilon_j)=0$ for any $1\leq i\leq m$. Using
(\ref{structureU}), this is equivalent to
$$\sum_{1\leq j,\ell \leq m}y_jc_{ij\ell}u^*(u_\ell)=0$$
for any $1\leq i\leq m$. This is the same with
\begin{equation}\label{sistemU}
\sum_{1\leq j\leq m}\left(\sum_{1\leq \ell \leq
m}u^*(u_l)c_{ij\ell}\right)y_j=0 \;\;\; \mbox{ for any }1\leq i\leq
m
\end{equation}
We conclude that $\varphi$ is injective if and only if the system
(\ref{sistemU}) has only the solution $y_1=\ldots=y_m=0$, which is
equivalent to the invertibility of the matrix $(\sum_{1\leq \ell\leq
m}u^*(u_\ell) c_{ij\ell})_{1\leq i,j\leq m}$. This proves the first
part of the proposition if we take $\alpha_\ell=u^*(u_\ell)$, and
also describes all the isomorphisms $\varphi$ as in the second part.
\end{proof}

In the rest of the paper we keep the assumptions from the
Introduction: $A$ is a finite-dimensional $G$-graded algebra,
$(e_j)_{j\in J}$ is a fixed homogeneous basis of $A$, and $\sigma$
is a fixed element of $G$. The grading on $A$ induces a $G$-graded
left $A$-module structure on $A^*$, whose homogeneous component of
degree $g$ is
$$(A^*)_g=\{ a^*\in A^*\,|\, a^*(A_h)=0 \mbox{ for any }h\neq
g^{-1}\}$$ for any $g\in G$. We also consider the shift $A(\sigma)$
of $A$, which is just the left $A$-module $A$, with grading defined
by $A(\sigma)_g=A_{g\sigma}$ for any $g\in G$. Then $A$ is called a
$\sigma$-graded Frobenius algebra if $A^*\simeq A(\sigma)$ as graded
left $A$-modules.\\

{\bf Proof of Proposition B:} $A$ is $\sigma$-graded Frobenius if
$A^*\simeq A(\sigma)$, or equivalently, if there exists an
isomorphism of graded left $A$-modules $\varphi:A\simeq
A^*(\sigma^{-1})$. Since in particular such a $\varphi$ is an
isomorphism of left $A$-modules from $A$ to $A^*$, we have by
Proposition \ref{propUdual} that $\varphi$ is of the form
$\varphi(a)=a\rightharpoonup a^*$ for some $a^*\in A^*$ such that
the matrix $(\sum_{\ell\in J}\alpha_\ell c_{ij\ell})_{i,j\in J}$ is
invertible, where $\alpha_\ell=a^*(e_\ell)$ for any $\ell\in J$.
Such a $\varphi$ is even an isomorphism of graded left $A$-modules
from $A$ to $A^*(\sigma^{-1})$ if and only if $a^*=\varphi(1)\in
A^*(\sigma^{-1})_\varepsilon=(A^*)_{\sigma^{-1}}$, which means that
$a^*(e_\ell)=0$ for any $\ell\in J\setminus J_\sigma$.

We conclude that $A$ is $\sigma$-graded Frobenius if and only if
there exists a family $(\alpha_\ell)_{\ell\in J_{\sigma}}\subset K$
such that the matrix $(\sum_{\ell\in J}\alpha_\ell
c_{ij\ell})_{i,j\in J}$ is invertible. This ends the proof, since
this matrix is just $P(\alpha,\sigma)$, where
$\alpha=(\alpha_\ell)_{\ell\in J_\sigma}$.   \qed \\

As for the graded symmetric property, we have the following.\\

{\bf Proposition B'.} {\it $A$ is graded symmetric if and only if
there exists a family of scalars $\alpha=(\alpha_\ell)_{\ell\in
J_\varepsilon}$ such that the $\sigma$-paratrophic matrix
$P(\varepsilon,\alpha)$ is symmetric and invertible.}\\[2mm]
{\bf Proof:} Keeping the notation from the proof of Proposition B
for $\sigma=\varepsilon$, an isomorphism $\varphi:A\ra A^*$ of
graded left $A$-modules is given by an $a^*\in (A^*)_\varepsilon$,
i.e., by a family of scalars $\alpha=(\alpha_\ell)_{\ell\in
J_\varepsilon}$, where $\alpha_\ell=a^*(e_\ell)$ for any $\ell\in
J_\varepsilon$, such that $P(\varepsilon,\alpha)$ is invertible.
Moreover, $\varphi$ is a morphism of $A$-bimodules if and only if
$e_i\rightharpoonup a^*=a^*\leftharpoonup e_i$ for any $i\in J$,
where $\leftharpoonup$ is the usual right $A$-action on $A^*$. This
is equivalent to $a^*(e_je_i)=a^*(e_ie_j)$ for any $i,j\in J$, which
is the same with $\sum_{\ell \in
J_\sigma}c_{ji\ell}\alpha_\ell=\sum_{\ell \in
J_\sigma}c_{ij\ell}\alpha_\ell$ for any $i,j$, i.e., the matrix
$P(\varepsilon,\alpha)$ is symmetric. \qed

\begin{remark}
{\rm If $A$ is a finite-dimensional algebra regarded with the
trivial $G$-grading, where $G$ is an arbitrary group, then $A$ is
$\varepsilon$-graded Frobenius (respectively graded symmetric) if
and only if it is a Frobenius algebra (respectively symmetric). In
this case, an $\varepsilon$-paratrophic matrix is just a usual
paratrophic matrix, and Proposition B and Proposition B' are just
the Frobenius' criterion and Nakayama-Nesbitt test for symmetric
algebras, see \cite[16.82 and Exercise 28, page 457]{lam2}. }
\end{remark}

\section{Koszul duals of quantum polynomial algebras and matrix algebras with good gradings} \label{sectiuneexemple}

We first recall that an algebra $A$ is a quadratic algebra if
$A\simeq T(V)/(R)$, where $T(V)$ is the tensor algebra of a
finite-dimensional space $V$, and $(R)$ is the ideal of $T(V)$
generated by a subspace $R$ of $T^2(V)=V\ot V$. The Koszul dual of
such an $A$ is the algebra $A^!=T(V^*)/(R^\perp)$, where $R^\perp$
is the orthogonal of $R$ in $V^*\ot V^*$.

If we take $A=K_{\bf q}[X_1,\ldots,X_n]$, the quantum polynomial
algebra associated with ${\bf q}=(q_{ij})_{1\leq i<j\leq n}$, then
$A\simeq T(V)/(R)$, where $V$ is a vector space with basis $\{
v_1,\ldots,v_n\}$, and $R$ is the subspace of $V\ot V$ spanned by
all $v_j\ot v_i-q_{ij}v_i\ot v_j$, with $1\leq i<j\leq n$.

Let $z=\sum_{1\leq k,\ell\leq n} \alpha_{k\ell}v_k^*\ot
v_{\ell}^*\in V^*\ot V^*$. Then for some $1\leq i<j\leq n$ we have
$z(v_j\ot v_i-q_{ij}v_i\ot v_j)=\alpha_{ji}-q_{ij}\alpha_{ij}$, so
$z\in R^\perp$ if and only if
$$z=\sum_{1\leq i<j\leq n}\alpha_{ij}(v_i^*\ot v_j^*+q_{ij}v_j^*\ot
v_i^*)+\sum_{1\leq i\leq n}\alpha_{ii}v_i^*\ot v_i^*$$ for some
$(\alpha_{ij})_{1\leq i\leq j\leq n}$. Thus $R^\perp$ is spanned by
all $v_i^*\ot v_j^*+q_{ij}v_j^*\ot v_i^*$ with $i<j$, and all
$v_i^*\ot v_i^*$. We conclude that
$$A^!\simeq K\langle X_1,\ldots,X_n\rangle /(X_jX_i+q_{ij}^{-1}X_iX_j,X_i^2|1\leq
i<j\leq n),$$ which can be rewritten as
$$K_{\bf q}[X_1,\ldots,X_n]^!\simeq A({\bf q'}), \mbox{  where }{\bf
q'}=(-q_{ij}^{-1})_{1\leq i<j\leq n}.$$

In the first part of this section we investigate Frobenius
properties of algebras of the type $A({\bf q})$, by proving
Proposition C. It is useful to label the basis elements of $A({\bf
q})$ by elements of $\mathbb{Z}_2^n$. Thus for any
$g=(\varepsilon_1,\ldots,\varepsilon_n)\in \mathbb{Z}_2^n$, denote
$e_g=x_1^{\varepsilon_1}\ldots x_n^{\varepsilon_n}$. For such a $g$,
we also denote ${\rm supp}(g)=\{ i|\varepsilon_i=\hat{1}\}$ and
$\overline{g}=(\hat{1}-\varepsilon_1,\ldots,\hat{1}-\varepsilon_n)$.

Let $g=(\varepsilon_1,\ldots,\varepsilon_n)$ and
$h=(\mu_1,\ldots,\mu_n)$ be elements of $\mathbb{Z}_2^n$. Then
$e_ge_h=0$ if there exists $1\leq i\leq n$ such that
$\varepsilon_i=\mu_i=\hat{1}$, i.e., if ${\rm supp}(g)\cap {\rm
supp}(h)\neq \emptyset$. In the case where ${\rm supp}(g)\cap {\rm
supp}(h)= \emptyset$, we have \bea
e_ge_h&=&(x_1^{\varepsilon_1}\ldots
x_n^{\varepsilon_n})(x_1^{\mu_1}\ldots x_n^{\mu_n})\\
&=& q_{12}^{\mu_1\varepsilon_2}\ldots
q_{1n}^{\mu_1\varepsilon_n}q_{23}^{\mu_2\varepsilon_3}\ldots
q_{2n}^{\mu_2\varepsilon_n}\ldots
q_{n-1,n}^{\mu_{n-1}\varepsilon_n}x_1^{\varepsilon_1+\mu_1}\ldots
x_n^{\varepsilon_n+\mu_n}\\
&=&\prod_{1\leq i<j\leq n}q_{ij}^{\mu_i\varepsilon_j}e_{g+h}\eea
Therefore, the structure constants associated with the basis
$(e_g)_{g\in \mathbb{Z}_2^n}$ are\\

$\bullet\;\;\; c_{g,h,\ell}=0$ for any $g,h,\ell$ such that ${\rm
supp}(g)\cap
{\rm supp}(h)\neq \emptyset$\\

$\bullet\;\;\; c_{g,h,\ell}=0$ for any $g,h,\ell$ such that ${\rm
supp}(g)\cap
{\rm supp}(h)= \emptyset$ and $\ell \neq g+h$\\

$\bullet\;\;\; c_{g,h,g+h}=\prod_{1\leq i<j\leq
n}q_{ij}^{\mu_i\varepsilon_j}$ if ${\rm
supp}(g)\cap {\rm supp}(h)= \emptyset$.\\

Now we prove Proposition C (ii). Thus $A({\bf q})$ is regarded as a
$\mathbb{Z}_2^n$-graded algebra whose homogeneous component of
degree $g=(\varepsilon_1,\ldots,\varepsilon_n)$ is
$Kx_1^{\varepsilon_1}\ldots x_n^{\varepsilon_n}$. We first show that
$A({\bf q})$ is $\tau$-graded Frobenius for
$\tau=(\hat{1},\ldots,\hat{1})$. Indeed, let $\alpha=(\alpha_\tau)$
with $\alpha_\tau\in K\setminus \{ 0\}$, and consider the
$\tau$-paratrophic matrix $P(\tau,\alpha)$, which by definition has
$c_{g,h,\tau}\alpha_\tau$ on position $(g,h)$ for any $g,h\in
\mathbb{Z}_2^n$. Taking into account the values of the structure
constants, for any $g\in \mathbb{Z}_2^n$, $P(\tau,\alpha)$ has
precisely one non-zero entry on row $g$, and that is on position
$(g,\overline{g})$. As $\overline{g}=\overline{h}$ if and only if
$g=h$, we see that $P(\tau,\alpha)$ is invertible, thus $A({\bf q})$
is $\tau$-graded Frobenius by Proposition B.

Now take $\sigma\neq (\hat{1},\ldots,\hat{1})$ and
$\alpha=(\alpha_\sigma)$ with $\alpha_\sigma\in K$, and let
$P(\sigma,\alpha)$ be the associated paratrophic matrix. Since
$\sigma$ has at least one component equal to $\hat{0}$, there exists
$g\in \mathbb{Z}_2^n$ such that ${\rm supp}(g)\cap {\rm
supp}(g+\sigma)\neq \emptyset$; indeed, one can pick a $g$ having
$\hat{1}$ on a position where $\sigma$ has $\hat{0}$. For such a
$g$, the $g$th row of $P(\sigma,\alpha)$ is all zero, since its
$(g,h)$-entry is $\alpha_\sigma c_{g,h,\sigma}$, which is 0 if
$g+h\neq \sigma$, while for $g+h=\sigma$, i.e., $h=g+\sigma$, we
have ${\rm supp}(g)\cap {\rm supp}(h)\neq \emptyset$, so again, the
$(g,h)$-entry is 0. We conclude that $P(\sigma,\alpha)$ cannot be
invertible, so $A({\bf
q})$ is not $\sigma$-graded Frobenius.\\

Next we prove Proposition C (iii). For any $0\leq k\leq n$ we have
$J_k=\{ g|\, |{\rm supp}(g)|=k\}$; in particular $J_n=\{ \tau\}$ and
$J_0=\{ 0\}$, where the element $(\hat{0},\ldots,\hat{0})$ of
$\mathbb{Z}_2^n$ is denoted by 0. Then $A({\bf q})$ is $n$-graded
Frobenius, with exactly the same proof as above to show that $A({\bf
q})$ is $\tau$-graded Frobenius.

Let now $0\leq k<n$. If $\alpha=(\alpha_\ell)_{\ell \in J_k}$ is a
family of scalars, then Theorem A shows that $P(k,\alpha)$ can not
be invertible if $k>0$, since $|J_k|={n\choose k}\neq 1=|J_0|$. If
$k=0$, the $(g,h)$-entry of $P(0,\alpha)$ is $\alpha_0c_{g,h,0}$,
which may be non-zero only for $g+h=0$, i.e., $g=h$, and ${\rm
supp}(g)\cap {\rm supp}(h)=\emptyset$, thus only for $g=h=0$. As it
has at most one non-zero entry, $P(0,\alpha)$ is not invertible,
thus $A({\bf q})$ is not 0-graded Frobenius. We conclude that
$A({\bf q})$ is not $k$-graded Frobenius for any $0\leq k<n$. It is
also not $k$-graded Frobenius for any $k<0$ or $k>n$, as
$P(k,\alpha)$ has only zero entries.\\

Finally, we prove Proposition C (i). As $A({\bf q})$ is
$\tau$-graded Frobenius, when regarded with the
$\mathbb{Z}_2^n$-grading, it is obviously Frobenius. In order to
investigate the symmetric property, we need to look closer at a
paratrophic matrix $P(\alpha)$, with $\alpha=(\alpha_\ell)_{\ell \in
\mathbb{Z}_2^n}$.

Assume first that $A({\bf q})$ is symmetric, thus there exists a
symmetric invertible matrix $P(\alpha)$. On row $\tau$ of
$P(\alpha)$ there may be a non-zero entry only on position
$(\tau,0)$, where we have the entry $\alpha_\tau$. As $P(\alpha)$ is
invertible, $\alpha_\tau$ must be non-zero. Then for any
$g=(\widehat{\varepsilon_1},\ldots,\widehat{\varepsilon_n})\in
\mathbb{Z}_2^n$, with $\varepsilon_1,\ldots,\varepsilon_n\in \{
0,1\}$, the $(g,\overline{g})$-entry of $P(\alpha)$ is
$$\alpha_\tau \prod_{1\leq i<j\leq
n}q_{ij}^{(1-\varepsilon_i)\varepsilon_j},$$ while the
$(\overline{g},g)$-entry is
$$\alpha_\tau
\prod_{1\leq i<j\leq n}q_{ij}^{\varepsilon_i(1-\varepsilon_j)}.$$
Since $P(\alpha)$ is symmetric and $\alpha_\tau\neq 0$, we get that
\begin{equation} \label{condsimetric0}
\prod_{1\leq i<j\leq
n}q_{ij}^{(1-\varepsilon_i)\varepsilon_j}=\prod_{1\leq i<j\leq
n}q_{ij}^{\varepsilon_i(1-\varepsilon_j)},
\end{equation}
or, equivalently,
\begin{equation}  \label{condsimetric1}
\prod_{1\leq i<j\leq n}q_{ij}^{\varepsilon_j-\varepsilon_i}=1.
\end{equation}
Thus (\ref{condsimetric1}) holds for any
$\varepsilon_1,\ldots,\varepsilon_n\in \{ 0,1\}$. If we fix an
$1\leq i< n$, and choose $\varepsilon_i=1$ and $\varepsilon_j=0$ for
any $j\neq i$, (\ref{condsimetric1}) becomes
$$\left(\prod_{1\leq p<i}q_{pi}\right)\cdot \left(\prod_{i<p\leq n}q_{ip}^{-1}\right)=1,$$
which is the same as
\begin{equation} \label{condsimetric2}
q_{in}=\left(\prod_{1\leq p<i}q_{pi}\right)\cdot
\left(\prod_{i<p\leq n-1}q_{ip}^{-1}\right).
\end{equation}

Conversely, assume that (\ref{condsimetric2}) holds for any $1\leq
i<n$. We show that $A({\bf q})$ is symmetric. If
$\varepsilon_1,\ldots,\varepsilon_n\in \{ 0,1\}$, then

\bea \prod_{1\leq i<j\leq
n}q_{ij}^{\varepsilon_j-\varepsilon_i}&=&\left(\prod_{1\leq i<j\leq
n-1}q_{ij}^{\varepsilon_j-\varepsilon_i}\right)\cdot
\left(\prod_{1\leq i\leq
n-1}q_{in}^{\varepsilon_n-\varepsilon_i}\right)\\
&=&\left(\prod_{1\leq i<j\leq
n-1}q_{ij}^{\varepsilon_j-\varepsilon_i}\right)\cdot
\left(\prod_{1\leq i\leq n-1 \atop 1\leq
p<i}q_{pi}^{\varepsilon_n-\varepsilon_i}\right)\cdot
\left(\prod_{1\leq i\leq n-1 \atop i<p \leq
n-1}q_{ip}^{\varepsilon_i-\varepsilon_n}\right)\\
&=&1. \eea The last equality holds since for each pair $(r,s)$ with
$1\leq r<s\leq n-1$, the exponents of $q_{rs}$ in the three factors
of the last product are $\varepsilon_s-\varepsilon_r$,
$\varepsilon_n-\varepsilon_s$ and $\varepsilon_r-\varepsilon_n$, and
their sum is zero. This shows that (\ref{condsimetric1}) holds for
any $\varepsilon_1,\ldots,\varepsilon_n\in \{ 0,1\}$. Now take the
family of scalars $\alpha=(\alpha_\ell)_{\ell\in \mathbb{Z}_2^n}$
such that $\alpha_\tau=1$ and $\alpha_\ell=0$ for any $\ell\neq
\tau$. Then for any $g,h\in \mathbb{Z}_2^n$, the $(g,h)$-entry of
the paratrophic matrix $P(\alpha)$ is $c_{g,h,\tau}$, which is
non-zero only if $g+h=\tau$, i.e., if $h=\overline{g}$. Then
obviously the $(g,h)$ and $(h,g)$ entries are equal for any $g,h$
such that $g+h\neq \tau$. On the other hand, the $(g,\overline{g})$
and $(\overline{g},g)$ entries of the paratrophic matrix $P(\alpha)$
are equal for any $g$, since equations (\ref{condsimetric1}) and
(\ref{condsimetric0}) are equivalent. Thus $P(\alpha)$ is symmetric,
and clearly invertible, and we conclude that $A({\bf q})$ is
symmetric.

\begin{remark}
{\rm Let $n\geq 2$ and let $V$ be a vector space of dimension $n$.
The exterior algebra $\wedge V$ is isomorphic to $A({\bf q})$, where
${\bf q}=(q_{ij})_{1\leq i<j\leq n}$ with $q_{ij}=-1$ for any $1\leq
i<j\leq n$. Proposition C shows that $\wedge V$ is a Frobenius
algebra, in fact it is even $n$-graded Frobenius when regarded with
the $\mathbb{Z}$-grading, and also $\tau$-graded Frobenius when
regarded as a $\mathbb{Z}_2^n$-graded algebra. Also, $\wedge V$ is
symmetric if and only $-1=(-1)^{n-2}$, i.e., if and only if $n$ is
odd.}
\end{remark}

At the end of this section we consider another class of graded
algebras: matrix algebras with good gradings. Let $A=M_n(K)$, where
$n$ is a positive integer. A good grading on $A$, by an arbitrary
group $G$, is a grading such that all matrix units $e_{ij}$, with
$1\leq i,j\leq n$, are homogeneous elements. It is showed in
\cite[Section 1]{dinr} that good $G$-gradings on $A$ are given by a
family $g_1,\ldots,g_n$ of elements of $G$, such that $e_{ij}$ has
degree $g_ig_j^{-1}$ for any $i,j$. For any $g\in G$, the
homogeneous component of degree $g$ of $A$ is $A_g=\langle e_{ij}|\,
g_ig_j^{-1}=g\rangle$. For such a grading, we consider the
homogeneous basis $(e_{ij})_{1\leq i,j\leq n}$, and keeping the
notation from the introduction, the associated sets of indices are
$J=\{(i,j)|\, 1\leq i,j\leq n\}$, and $J_g=\{ (i,j)|\,
g_ig_j^{-1}=g\}$ for any $g\in G$; in particular, $J_\varepsilon=\{
(i,j)|\, g_i=g_j\}$. Since $e_{ij}e_{pq}=\delta_{jp}e_{iq}$ for any
$i,j,p,q$, the structure constants associated with this basis are
$c_{(i,j),(p,q),(r,s)}=1$ if $j=p$ and $(r,s)=(i,q)$, and
$c_{(i,j),(p,q),(r,s)}=1$ in any other case.

Let us consider the family $\alpha=(\alpha_{k\ell})_{(k,\ell)\in
J_\varepsilon}$ such that $\alpha_{kk}=1$ for any $1\leq k\leq n$,
and $\alpha_{k\ell}=0$ for any $(k,\ell)\in J_\varepsilon$ with
$k\neq \ell$. The $((i,j),(p,q))$-entry of the
$\varepsilon$-paratrophic matrix $P(\varepsilon,\alpha)$ is
$$\sum_{(k,\ell)\in
J_\varepsilon}\alpha_{k\ell}c_{(i,j),(p,q),(k,\ell)} =\left\{
\begin{array}{ll}
1,&\mbox{ if } (p,q)=(j,i),\\
0,&\mbox{ otherwise. }
\end{array}
\right.$$ Therefore, $P(\varepsilon,\alpha)$ has precisely one
non-zero entry on each row and on each column, thus it is
invertible, and moreover, it is symmetric. This shows that $A$ is
graded symmetric, a fact that was remarked in \cite[Example
6.4]{dnn1}, using a different approach.

We note that good $G$-gradings on $M_n(K)$ are classified in
\cite[Corollary 2.2]{cdn} by the orbits of the biaction of the
symmetric group $\mathcal{S}_n$ by permutation, and of $G$ by right
translation, on the set $G^n$.

\section{A structure result for $\sigma$-graded Frobenius
algebras}\label{sectionstructureresult}

We recall that the $G$-graded algebra is called left
$\sigma$-faithful if for any $g\in G$ and any $a\in A_g$ such that
$A_{\sigma g^{-1}}a=0$, we must have $a=0$.

\begin{proposition}  \label{propfaithful}
$A$ is left $\sigma$-faithful if and only if ${\rm rank}\,
C_g=|J_g|$ for any $g\in G$.
\end{proposition}
\begin{proof}
Let $g\in G$ and $a=\sum_{j\in J_g}y_je_j\in A_g$. For any $r\in
J_{\sigma g^{-1}}$ we have \bea e_ra&=&\sum_{j\in J_g}y_je_re_j\\
&=&\sum_{j\in J_g}\sum_{\ell \in J_\sigma}y_jc_{rj\ell}e_\ell\\
&=&\sum_{\ell\in J_\sigma}\left(\sum_{j\in
J_g}y_jc_{rj\ell}\right)e_\ell \eea Then $A_{\sigma g^{-1}}a=0$ if
and only if $e_ra=0$ for any $r\in J_{\sigma g^{-1}}$, and by the
above computation, this is equivalent to $\sum_{j\in
J_g}y_jc_{rj\ell}=0$ for any $r\in J_{\sigma g^{-1}}$ and $\ell \in
J_\sigma$. But this is a linear system in $(y_j)_{j\in J_g}$ whose
matrix is $C_g$, thus it has only the trivial solution if and only
if ${\rm rank}\, C_g=|J_g|$. We conclude that the only $a\in A_g$
such that $A_{\sigma g^{-1}}a=0$ is $a=0$ if and only if ${\rm
rank}\, C_g=|J_g|$, and this ends the proof.
\end{proof}

As a particular case of Proposition \ref{propUdual} we have the
following.

\begin{corollary} \label{corAsigma}
$(A_\sigma)^*$ is isomorphic to $A_\varepsilon$ as a left
$A_\varepsilon$-module if and only if $|J_\sigma|=|J_\varepsilon|$
and there exists a family $(\alpha_\ell)_{\ell \in J_\sigma}\subset
K$ such that the matrix $(\sum_{\ell\in J_\sigma}\alpha_\ell
c_{ij\ell})_{i\in J_\sigma, j\in J_\varepsilon}$ is invertible.
\end{corollary}

{\bf Proof of Corollary D:} Assume that $A$ is graded
$\sigma$-Frobenius. By Proposition B, there exists a family
$\alpha=(\alpha_\ell)_{\ell\in J_\sigma}\subset K$ such that
$P(\alpha,\sigma)$ is invertible. Using Theorem A, we have that
$|J_\varepsilon|=|J_\sigma|$ and
$P(\sigma,\alpha)_{J_\varepsilon,J_\sigma}$ is invertible, so
$(A_\sigma)^*\simeq A_\varepsilon$ as left $A_\varepsilon$-modules
by Corollary \ref{corAsigma}, and also ${\rm rank}\, C_g=|J_g|$ for
any $g\in G$, so $A$ is left $\sigma$-faithful by Proposition
\ref{propfaithful}.

Conversely, assume that $A$ is left $\sigma$-faithful and
$(A_\sigma)^*\simeq A_\varepsilon$ as left $A_\varepsilon$-modules.
By Corollary \ref{corAsigma} we have $|J_\sigma|=|J_\varepsilon|$
and there exists a family $(\alpha_\ell)_{\ell \in J_\sigma}\subset
K$ such that the matrix $(\sum_{\ell\in J_\sigma}\alpha_\ell
c_{ij\ell})_{i\in J_\sigma, j\in J_\varepsilon}$ is invertible.
Denote $\alpha=(\alpha_\ell)_{\ell \in J_\sigma}$. Then the latter
mentioned matrix is just
$P(\sigma,\alpha)_{J_\varepsilon,J_\sigma}$, and it is invertible.
Using Proposition \ref{propfaithful}, ${\rm rank}\, C_g=|J_g|$ for
any $g\in G$. Now Theorem A shows that $P(\sigma,\alpha)$ is
invertible, and then $A$ is $\sigma$-graded Frobenius by Proposition
B. \qed \\

{\bf Acknowledgement.} The first and the fourth authors were
supported by the PNRR grant CF 44/14.11.2022 {\it Cohomological Hall
algebras of smooth surfaces and applications}.


\begin{thebibliography}{99}

\bibitem{cdn}
S. Caenepeel, S. D\u{a}sc\u{a}lescu, C. N\u{a}st\u{a}sescu, On
gradings of matrix algebras and descent theory, Comm. Algebra {\bf
30} (2002), 5901-5920.

\bibitem{dinr}
S. D\u{a}sc\u{a}lescu, B. Ion, C. N\u{a}st\u{a}sescu, J. Rios
Montes,  Group gradings on full matrix rings, J. Algebra {\bf 220}
(1999), 709-728.


\bibitem{dnn1} S. D\u{a}sc\u{a}lescu, C. N\u{a}st\u{a}sescu and L.
N\u{a}st\u{a}sescu, Frobenius algebras of corepresentations and
group graded vector spaces, J. Algebra {\bf 406} (2014), 226-250.


\bibitem{Frobenius} F. G. Frobenius, Theorie der hyperkomplexen Gr${\rm\ddot{o}}\beta$en I. II., Berl. Ber. 1903 (1903), 504-537, 634-645.



\bibitem{jor} P. J$\rm \o$rgensen, A noncommutative BGG correspondence, Pacific J. Math. {\bf 218} (2005), 357-377.




\bibitem{lam2} T. Y. Lam, Lectures on modules and rings, GTM {\bf
189}, Springer Verlag, 1999.

\bibitem{lpwz} D.-M. Lu, J. H. Palmieri, Q.-S. Wu and J. J. Zhang,  Regular algebras of dimension
4 and their $A_\infty$-Ext-algebras, Duke Math. J. {\bf 137} (2007),
537-584.

\bibitem{nvo} C. N\u{a}st\u{a}sescu and F. van Oystaeyen, Methods
of graded rings, Lecture Notes in Math., vol. 1836 (2004), Springer
Verlag.




\end{thebibliography}
\end{document}